\documentclass[11pt]{amsart}
\usepackage{amsmath,amssymb,amsthm}
\usepackage{graphicx}
\usepackage{booktabs}
\usepackage{array}
\usepackage[margin=1in]{geometry}
\usepackage{xcolor}
\usepackage[colorlinks=true,linkcolor=blue!50!black,citecolor=blue!50!black,urlcolor=blue!50!black]{hyperref}

\newtheorem{theorem}{Theorem}
\newtheorem{proposition}[theorem]{Proposition}
\newtheorem{lemma}[theorem]{Lemma}
\newtheorem{corollary}[theorem]{Corollary}
\theoremstyle{definition}
\newtheorem{definition}[theorem]{Definition}
\theoremstyle{remark}
\newtheorem{remark}[theorem]{Remark}
\newtheorem{example}[theorem]{Example}

\newcommand{\Sodd}{\Sigma_{\mathrm{odd}}}
\newcommand{\Seven}{\Sigma_{\mathrm{even}}}
\DeclareMathOperator{\Imag}{Im}

\begin{document}

\title[Equal area sums of alternate sectors]{Equating the area sums of alternate sectors of a disc cut by concurrent chords}

\author {Azhar Iqbal and Derek Abbott}

\address{School of Electrical and Mechanical Engineering, College of Engineering and Information Technology, Adelaide University,
Adelaide, South Australia 5005, Australia}

\subjclass[2020]{51M25, 51M04, 26B15}
\keywords{Pizza theorem, curvilinear sectors of a disc, concurrent chords, power of a point, fair division}

\begin{abstract}
A disc is cut by $n$ chords through a common interior point $O$, which produces $2n$ curvilinear sectors. We determine when the two sums of alternate sectors have equal area for arbitrary chord directions, which need not be equally spaced. We first derive a closed form for the area of a single sector. We then show that the answer depends on the parity of $n$. For even $n$, at fixed orientation the condition is affine in $r_0^2$, where $r_0$ is the distance from $O$ to the centre; it involves only the total angle of the odd sectors and one trigonometric moment of the directions. For odd $n$, the condition is a nonlinear alternating sum of circular-segment functions of the distances of the chords from the centre. We treat four, six and eight pieces in detail. We also hold the chord directions fixed and translate the whole pencil to an arbitrary interior point, and we characterise the direction sets that give equal sums at every such point. For four chords only equally spaced directions have this property, which gives a converse of the pizza theorem; for an odd number of chords no direction set has it. A fair two-chord cut requires the acute angle between the chords to lie within $(r_0/a)^2$ radians of a right angle, where $a$ is the radius. For equally spaced chords the results reduce to the known theorems of Goldberg and of Mabry and Deiermann.
\end{abstract}

\maketitle

\section{Introduction}

Problems in circle geometry~\cite{Ogilvy,AltshillerCourt,CoxeterGreitzer} often have answers that are neither intuitive nor obvious. A celebrated example is the \emph{pizza theorem}. Upton~\cite{Upton} proposed, and Goldberg~\cite{Goldberg} solved, the following problem. A disc is cut by four chords through an arbitrary interior point $O$, with consecutive chords meeting at angles of $\pi/4$; then the eight resulting pieces, coloured alternately, give two sets of equal total area. Goldberg also showed that the same holds for any even number $n\ge 4$ of equally spaced chords. The result has since received a proof without words~\cite{CarterWagon}, elementary analytic proofs~\cite{Hirschhorn,Pearce,Shirali}, popular expositions~\cite{Berzsenyi1,Berzsenyi2,KVW} and further dissection proofs~\cite{Frederickson}. For an odd number of equally spaced chords the two sums are in general unequal. The question of which share is larger was raised in~\cite{CarterDuncanWagon,DeiermannMabry} and settled by Mabry and Deiermann~\cite{MabryDeiermann}. Gluck~\cite{Gluck} later obtained an exact Fourier series for the imbalance.

\subsection*{Other generalisations, and how this paper differs}
The higher-dimensional pizza theorems of Brailov~\cite{Brailov} and of Ehrenborg, Morel and Readdy~\cite{EMR2022,EMR2023} (see also~\cite{Ehrenborg2025}) work with \emph{Coxeter} arrangements, that is, hyperplane arrangements that are invariant under reflection in each of their hyperplanes. In exchange, they allow any body with the same reflection symmetry in place of the disc. In the plane, a pencil of lines through a point that is invariant under reflection in each of its lines is necessarily equally spaced. Berele and Catoiu~\cite{BereleCatoiu} studied a polygonal (non-circular) analogue. In all of these works the cut arrangement is highly symmetric. Here the region is always a disc, but the pencil of chords is \emph{arbitrary}: the directions are any $\theta_1<\theta_2<\dots<\theta_n<\theta_1+\pi$. The reflection and dissection arguments of the symmetric theory are then unavailable, and the answer genuinely depends on where $O$ lies. We derive necessary and sufficient equations that determine which directions, positions of $O$ and radii make the alternate area sums equal. These equations characterise equality completely, although we classify their solutions explicitly only in particular cases (for example, $n=2$, $n=3$ with a chord through the centre, and universal fairness for $n=4$ and for odd $n$).

\subsection*{What is new}
Table~\ref{tab:novelty} lists each result and its status relative to the literature cited above. We use ``new'' to mean that we did not find the statement in that literature; our search was targeted and does not claim to be exhaustive. In brief, the main contributions are as follows.
\begin{enumerate}
\item Necessary and sufficient conditions for equal alternate sums with \emph{arbitrary} chord directions, whose form depends on the parity of $n$ (Theorems~\ref{thm:even} and~\ref{thm:odd}). These rest on two elementary identities that are known for equally spaced cuts and that we show to hold for arbitrary directions (Lemmas~\ref{lem:pair} and~\ref{lem:half}).
\item A criterion for \emph{universal fairness}, in the sense of Definition~\ref{def:universal} (Proposition~\ref{prop:universal}). From it we obtain a rigidity theorem for four chords (Proposition~\ref{prop:n4unique}) and the impossibility of universal fairness for any odd number of chords (Proposition~\ref{prop:oddnever}).
\item A detailed treatment of four, six and eight pieces. This includes fairness along a diameter (Section~\ref{sec:eight}), the admissible angle window for a fair two-chord cut (Proposition~\ref{prop:n2window}), and the three-chord conditions (Corollary~\ref{cor:n3}).
\end{enumerate}
Appendix~\ref{app:numerics} documents reproducible numerical checks of the principal area identities and of selected special cases.

\begin{table}[t]
\centering
\small
\renewcommand{\arraystretch}{1.15}
\begin{tabular}{>{\raggedright}p{0.38\textwidth}>{\raggedright\arraybackslash}p{0.54\textwidth}}
\toprule
Result & Status \\
\midrule
Prop.~\ref{prop:sector} (area of one sector) & Elementary calculus; given in closed form for reference. \\
Lemma~\ref{lem:pair} (opposite pair) & Arbitrary-angle form of the paired-sector identity behind~\cite{Goldberg,Hirschhorn}. \\
Lemma~\ref{lem:half}, Thm.~\ref{thm:odd} (odd $n$) & Arbitrary-angle form of the alternating strip identity of~\cite{MabryDeiermann}, which is stated there for equally spaced cuts. \\
Thm.~\ref{thm:even} (even $n$, arbitrary angles) & New as stated; follows directly from Lemma~\ref{lem:pair}. \\
Fairness along a diameter (Sec.~\ref{sec:eight}); Prop.~\ref{prop:universal} & New. \\
Cor.~\ref{cor:goldberg} & Known (Goldberg~\cite{Goldberg}). \\
Prop.~\ref{prop:n4unique} (rigidity for $n=4$) & New. \\
Cor.~\ref{cor:perp}, Prop.~\ref{prop:n2window} (two chords) & New; Cor.~\ref{cor:perp} is elementary. \\
Cor.~\ref{cor:n3} (three chords, arbitrary angles) & New. \\
Prop.~\ref{prop:asym} (leading-order imbalance, odd $n$) & Leading term of Gluck's exact series~\cite{Gluck}; the global sign rule is due to~\cite{MabryDeiermann}. We add a short derivation with an explicit, uniform remainder bound. \\
Prop.~\ref{prop:oddnever} (no universally fair odd cut) & For equally spaced cuts this is implied by~\cite{MabryDeiermann}; new for arbitrary directions. \\
\bottomrule
\end{tabular}
\caption{The results of this paper and their relation to the literature.}
\label{tab:novelty}
\end{table}

\section{Set-up and the area of a single sector}\label{sec:setup}

Place the pole of a polar coordinate system at the interior point $O$, and let the centre $C$ of the circle of radius $a$ have polar coordinates $(r_0,\theta_0)$, with $0\le r_0<a$. We write
\[
\kappa=\frac{r_0}{a}\in[0,1),\qquad \varphi=\theta-\theta_0 .
\]
The ray from $O$ in direction $\theta$ meets the circle at distance
\begin{equation}\label{eq:r}
r(\theta)=r_0\cos(\theta-\theta_0)+\sqrt{a^{2}-r_0^{2}\sin^{2}(\theta-\theta_0)} .
\end{equation}
This follows from the law of cosines in triangle $OCP$, where $P$ is the point on the circle. The square root is real and $r(\theta)>0$ because $r_0<a$. Two identities follow immediately from~\eqref{eq:r}:
\begin{equation}\label{eq:power}
r(\theta)\,r(\theta+\pi)=a^{2}-r_0^{2},\qquad r(\theta)-r(\theta+\pi)=2r_0\cos(\theta-\theta_0).
\end{equation}
The first is the intersecting-chords (power of a point) theorem~\cite[\S2.1]{CoxeterGreitzer}.

\begin{figure}[t]
\centering
\includegraphics[width=\textwidth]{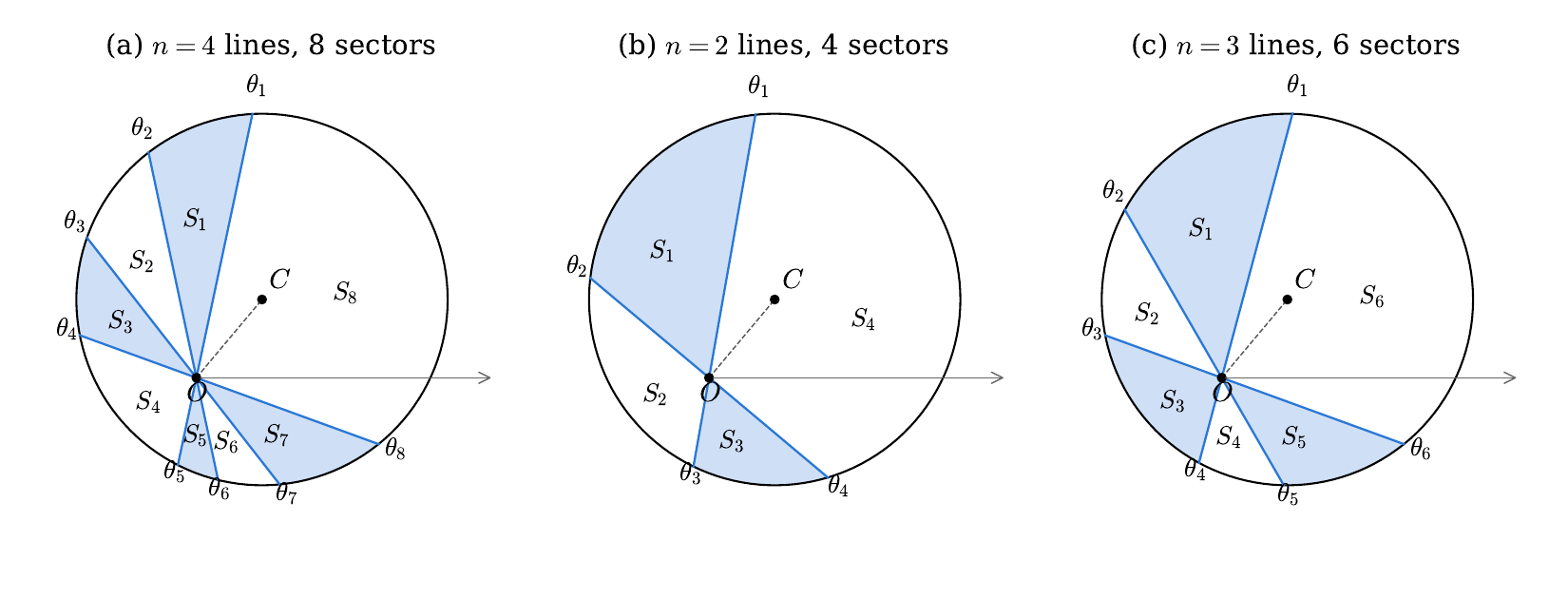}
\caption{A disc with centre $C=(r_0,\theta_0)$ cut by $n$ chords through the pole $O$. The odd-numbered pieces are shaded. (a)~Four chords (eight pieces); (b)~two chords (four pieces); (c)~three chords (six pieces). In all panels $r_0/a=0.55$ and $\theta_0=50^\circ$.}
\label{fig:geometry}
\end{figure}

Let $n\ge1$ chords pass through $O$ in directions $\theta_1<\theta_2<\dots<\theta_n<\theta_1+\pi$. Extend the list by $\theta_{j+n}=\theta_j+\pi$ for $1\le j\le n$ and $\theta_{2n+1}=\theta_1+2\pi$. The chords divide the disc into $2n$ \emph{curvilinear sectors} (pieces)
\[
S_j=\{\text{points of the disc with polar angle in }[\theta_j,\theta_{j+1}]\},\qquad 1\le j\le 2n,
\]
numbered counter-clockwise (Figure~\ref{fig:geometry}). The apex of each piece is $O$, which in general is not the centre of the circle. We use $S_j$ both for the piece and for its area, write $\varphi_j=\theta_j-\theta_0$, and set
\[
\Sodd=S_1+S_3+\dots+S_{2n-1},\qquad \Seven=S_2+S_4+\dots+S_{2n},\qquad \Sodd+\Seven=\pi a^{2}.
\]
The two sums are therefore equal precisely when $\Sodd=\pi a^{2}/2$.

\begin{proposition}[Area of a sector]\label{prop:sector}
Let $h(s)=\arcsin s+s\sqrt{1-s^{2}}$ for $|s|\le1$, and define
\begin{equation}\label{eq:F}
F(\varphi)=\frac{a^{2}}{2}\Bigl[\varphi+h(\kappa\sin\varphi)\Bigr]+\frac{r_0^{2}}{4}\sin2\varphi .
\end{equation}
Then the area of the region of the disc swept by the rays from $O$ with polar angles between $\alpha$ and $\beta$ ($\alpha<\beta$) is
\begin{equation}\label{eq:sectorarea}
\frac12\int_{\alpha}^{\beta}r^{2}(\theta)\,d\theta=F(\beta-\theta_0)-F(\alpha-\theta_0).
\end{equation}
In particular, $S_j=F(\varphi_{j+1})-F(\varphi_j)$.
\end{proposition}

\begin{proof}
\emph{Case $\kappa=0$.} Then $r(\theta)=a$, the left side of~\eqref{eq:sectorarea} is $\tfrac12a^2(\beta-\alpha)$, and $F(\varphi)=\tfrac12a^2\varphi$ because $h(0)=0$ and $r_0=0$. The claim holds.

\emph{Case $0<\kappa<1$.} Squaring~\eqref{eq:r} gives
\[
\tfrac12 r^{2}=\tfrac12 a^{2}+\tfrac12 r_0^{2}\cos2\varphi+ar_0\cos\varphi\sqrt{1-\kappa^{2}\sin^{2}\varphi}.
\]
The first two terms integrate to $\tfrac12a^{2}\varphi+\tfrac14r_0^{2}\sin2\varphi$. For the last term, substitute $\sin x=\kappa\sin\varphi$. Because $|\kappa\sin\varphi|<1$, we have $x\in(-\pi/2,\pi/2)$, $\cos x=\sqrt{1-\kappa^{2}\sin^{2}\varphi}>0$ and $\cos x\,dx=\kappa\cos\varphi\,d\varphi$. Since $\kappa\neq0$,
\[
ar_0\int\cos\varphi\sqrt{1-\kappa^{2}\sin^{2}\varphi}\,d\varphi=\frac{ar_0}{\kappa}\int\cos^{2}x\,dx=\frac{a^{2}}{2}\Bigl(x+\tfrac12\sin2x\Bigr)=\frac{a^{2}}{2}\,h(\kappa\sin\varphi),
\]
where we used $\tfrac12\sin2x=\sin x\cos x=\kappa\sin\varphi\sqrt{1-\kappa^2\sin^2\varphi}$. Differentiating~\eqref{eq:F} directly confirms that $F'(\varphi)=\tfrac12r^{2}$ for every $\varphi$, using $h'(s)=2\sqrt{1-s^2}$. So $F$ is a global antiderivative, and no case distinction according to the sign of $\cos\varphi$ is needed.
\end{proof}

\begin{remark}\label{rem:coefficient}
The term in $x=\arcsin(\kappa\sin\varphi)$ carries the coefficient $a^{2}/2$. A coefficient involving $a/r_0$ would be dimensionally inconsistent and would diverge as $r_0\to0$.
\end{remark}

Two consequences of Proposition~\ref{prop:sector} drive everything that follows. Since $\sin2(\varphi+\pi)=\sin2\varphi$ and $h(\kappa\sin(\varphi+\pi))=-h(\kappa\sin\varphi)$ ($h$ is odd), we have
\begin{equation}\label{eq:Fshift}
F(\varphi+\pi)+F(\varphi)=a^{2}\varphi+\frac{\pi a^{2}}{2}+\frac{r_0^{2}}{2}\sin2\varphi,
\qquad
F(\varphi+\pi)-F(\varphi)=\frac{\pi a^{2}}{2}-a^{2}h(\kappa\sin\varphi).
\end{equation}

\begin{lemma}[Opposite sectors]\label{lem:pair}
For $1\le j\le n$, the vertically opposite pieces $S_j$ and $S_{j+n}$ satisfy
\begin{equation}\label{eq:pair}
S_j+S_{j+n}=a^{2}(\theta_{j+1}-\theta_j)+\frac{r_0^{2}}{2}\bigl[\sin2\varphi_{j+1}-\sin2\varphi_j\bigr].
\end{equation}
\end{lemma}

\begin{proof}
By Proposition~\ref{prop:sector}, $S_j+S_{j+n}=[F(\varphi_{j+1})+F(\varphi_{j+1}+\pi)]-[F(\varphi_j)+F(\varphi_j+\pi)]$. Apply the first identity in~\eqref{eq:Fshift}. For $j=n$ we have $\varphi_{n+1}=\varphi_1+\pi$, and~\eqref{eq:pair} still holds because $\sin2\varphi$ has period $\pi$.
\end{proof}

\begin{remark}
Lemma~\ref{lem:pair} also has a proof that avoids $\arcsin$ entirely. By~\eqref{eq:power},
$r^{2}(\theta)+r^{2}(\theta+\pi)=[r(\theta)-r(\theta+\pi)]^{2}+2r(\theta)r(\theta+\pi)=2a^2+2r_0^2\cos2\varphi$.
Integrating $\tfrac12(\cdot)$ from $\theta_j$ to $\theta_{j+1}$ gives~\eqref{eq:pair}. This is closely related to the Pythagorean argument of Hirschhorn et al.~\cite{Hirschhorn} for equally spaced chords.
\end{remark}

\begin{lemma}[Half-discs]\label{lem:half}
Let $H_j=S_j\cup S_{j+1}\cup\dots\cup S_{j+n-1}$ be the part of the disc to the left of the $j$th chord directed along $\theta_j$ (polar angles in $[\theta_j,\theta_j+\pi]$). Define the \emph{offset} $d_j=r_0\sin\varphi_j$. Then $|d_j|$ is the distance from $C$ to the $j$th chord, and $d_j>0$ exactly when $C$ lies strictly to the \emph{right} of the directed chord, that is, outside $H_j$. (This is the negative of the usual signed distance measured along the left normal.) With this convention,
\begin{equation}\label{eq:half}
|H_j|=\frac{\pi a^{2}}{2}-a^{2}h\Bigl(\frac{d_j}{a}\Bigr)=\frac{\pi a^{2}}{2}-a^{2}h(\kappa\sin\varphi_j),
\end{equation}
and $a^{2}h(|d_j|/a)$ is half the difference between the areas of the two pieces into which the $j$th chord cuts the disc.
\end{lemma}

\begin{proof}
Formula~\eqref{eq:half} is the second identity in~\eqref{eq:Fshift}. The vector from $O$ to $C$ is $r_0(\cos\theta_0,\sin\theta_0)$. Its component along the left unit normal $(-\sin\theta_j,\cos\theta_j)$ of the directed chord is $r_0\sin(\theta_0-\theta_j)=-d_j$, which gives the sign convention. A chord at distance $d\ge0$ from the centre cuts off a segment of area $a^{2}\arccos(d/a)-d\sqrt{a^{2}-d^{2}}$. Since $\arccos s=\pi/2-\arcsin s$, half the difference between the two pieces is $a^{2}h(d/a)$. The smaller piece is the one not containing $C$, consistent with the sign in~\eqref{eq:half}.
\end{proof}

\section{The general conditions: even versus odd numbers of chords}\label{sec:general}

\begin{theorem}[Even number of chords]\label{thm:even}
Let $n\ge2$ be even, and let $A=\sum_{j\ \mathrm{odd},\,j<n}(\theta_{j+1}-\theta_j)$ be the total opening angle of the odd pieces $S_1,S_3,\dots,S_{n-1}$. Then
\begin{equation}\label{eq:even}
\Sodd=a^{2}A+\frac{r_0^{2}}{2}\sum_{j=1}^{n}(-1)^{j}\sin2\varphi_j,
\qquad
\Sodd-\Seven=2a^{2}\Bigl(A-\frac{\pi}{2}\Bigr)+r_0^{2}\sum_{j=1}^{n}(-1)^{j}\sin2\varphi_j .
\end{equation}
Hence $\Sodd=\Seven$ if and only if
\begin{equation}\label{eq:evencond}
a^{2}\Bigl(A-\frac{\pi}{2}\Bigr)+\frac{r_0^{2}}{2}\sum_{j=1}^{n}(-1)^{j}\sin2(\theta_j-\theta_0)=0 .
\end{equation}
Equivalently, with $Z=\sum_{j=1}^{n}(-1)^{j}e^{2i\theta_j}$, the condition reads $a^{2}(A-\pi/2)+\tfrac12r_0^{2}\Imag\bigl(e^{-2i\theta_0}Z\bigr)=0$.
\end{theorem}

\begin{proof}
Because $n$ is even, $j$ and $j+n$ have the same parity. So the odd pieces come in opposite pairs $\{S_j,S_{j+n}\}$ with $j$ odd and $j\le n-1$. Summing~\eqref{eq:pair} over these $j$ gives the first formula in~\eqref{eq:even}: each odd index contributes $-\sin2\varphi_j$ and each even index $+\sin2\varphi_j$. The second formula follows from $\Seven=\pi a^{2}-\Sodd$.
\end{proof}

\begin{theorem}[Odd number of chords]\label{thm:odd}
Let $n\ge1$ be odd. Then
\begin{equation}\label{eq:odd}
\Sodd=\frac{\pi a^{2}}{2}-a^{2}\sum_{j=1}^{n}(-1)^{j+1}h(\kappa\sin\varphi_j),
\qquad
\Sodd-\Seven=-2a^{2}\sum_{j=1}^{n}(-1)^{j+1}h(\kappa\sin\varphi_j).
\end{equation}
Hence $\Sodd=\Seven$ if and only if
\begin{equation}\label{eq:oddcond}
\sum_{j=1}^{n}(-1)^{j+1}h\Bigl(\frac{r_0}{a}\sin(\theta_j-\theta_0)\Bigr)=0 .
\end{equation}
\end{theorem}

\begin{proof}
For $n$ odd, the indicator function of the odd pieces equals $\sum_{j=1}^{n}(-1)^{j+1}\mathbf 1_{H_j}$ away from the chords. To see this, take a point in $S_m$. If $1\le m\le n$, it lies in $H_j$ exactly when $j\le m$, and $\sum_{j\le m}(-1)^{j+1}$ is $1$ for $m$ odd and $0$ for $m$ even. If $m=n+l$ with $1\le l\le n$, it lies in $H_j$ exactly when $j\ge l+1$. Since $n$ is odd, $\sum_{j=l+1}^{n}(-1)^{j+1}$ is $1$ when $l$ is even (so $m$ is odd) and $0$ when $l$ is odd (so $m$ is even). Integrating and using~\eqref{eq:half} gives $\Sodd=\sum_j(-1)^{j+1}|H_j|$. This sum has $(n+1)/2$ plus signs and $(n-1)/2$ minus signs, so the constant parts contribute $\pi a^{2}/2$, which yields~\eqref{eq:odd}.
\end{proof}

\begin{remark}[Why the parity matters]\label{rem:parity}
For $n$ even, $\Sodd$ is a sum of opposite pairs, and the transcendental terms $h(\kappa\sin\varphi_j)$ cancel within each pair. At a fixed orientation the condition~\eqref{eq:evencond} is therefore \emph{affine in $r_0^2$}. For $n$ odd, every opposite pair contains one odd and one even piece. The elementary terms then cancel: the odd pieces always have total opening angle $\pi$, and the $\sin2\varphi$ terms telescope to zero. What remains is the \emph{nonlinear}, odd, real-analytic function of $\kappa=r_0/a$ in~\eqref{eq:oddcond}. That function depends only on the offsets $d_j$ of the chords from the centre. For equally spaced chords, $\sum_j(-1)^{j+1}h(d_j/a)$ is the alternating ``strip'' sum used by Mabry and Deiermann~\cite{MabryDeiermann}.
\end{remark}

\begin{remark}[Centred case]\label{rem:centred}
When $r_0=0$, the condition becomes $A=\pi/2$ for $n$ even and holds \emph{identically} for $n$ odd, because each sum is then the union of $n$ circular sectors of total angle $\pi$. In particular, for an odd number of chords through the centre the alternate sums are always equal, whatever the angles.
\end{remark}

\section{Eight pieces ($n=4$) and universal fairness}\label{sec:eight}

For four chords, Theorem~\ref{thm:even} gives
\begin{equation}\label{eq:n4}
\Sodd-\Seven=2a^{2}\Bigl(\theta_2-\theta_1+\theta_4-\theta_3-\frac{\pi}{2}\Bigr)+r_0^{2}\bigl[\sin2\varphi_2-\sin2\varphi_1+\sin2\varphi_4-\sin2\varphi_3\bigr],
\end{equation}
so the alternate sums are equal if and only if
\begin{equation}\label{eq:n4cond}
\frac{r_0^{2}}{2}\bigl[\sin2\varphi_2-\sin2\varphi_1+\sin2\varphi_4-\sin2\varphi_3\bigr]+a^{2}\Bigl(\theta_2-\theta_1+\theta_4-\theta_3-\frac{\pi}{2}\Bigr)=0 .
\end{equation}
When $O=C$ this reduces to $(\theta_2-\theta_1)+(\theta_4-\theta_3)=\pi/2$.

\subsection*{Moving the point of concurrency}
From now on we will often vary the point $O$ while keeping the chord \emph{directions} fixed. In other words, the pencil of lines is translated rigidly so that it passes through the new point, and the lines are not held fixed in the plane. (Fixed lines through a common point cannot pass through a different point.) The disc, and hence $a$ and $C$, stay fixed throughout. As $O$ ranges over the interior of the disc, the pair $(r_0,\theta_0)$ ranges over $[0,a)\times\mathbb R$, and the labels and parities of the pieces do not change.

\begin{definition}\label{def:universal}
A set of directions $\theta_1<\dots<\theta_n<\theta_1+\pi$ is \emph{universally fair} if, for \emph{every} interior point $O$ of the disc, the $n$ chords through $O$ in these directions give $\Sodd=\Seven$.
\end{definition}

\subsection*{Fairness along a diameter}
Fix the directions and $C$, and move $O$ along the diameter through $C$ in direction $\theta_0$. Both conditions below are unchanged under $\theta_0\mapsto\theta_0+\pi$, which covers points on either side of $C$. By~\eqref{eq:n4cond}, the sums are equal for \emph{every} point $O$ on that diameter, that is, for every $r_0\in[0,a)$, if and only if
\begin{equation}\label{eq:n4special}
(\theta_2-\theta_1)+(\theta_4-\theta_3)=\frac{\pi}{2}
\quad\text{and}\quad
\sin2\varphi_2+\sin2\varphi_4=\sin2\varphi_1+\sin2\varphi_3 .
\end{equation}
Using $\sin u+\sin v=2\sin\frac{u+v}{2}\cos\frac{u-v}{2}$ and substituting $\theta_2+\theta_4=\theta_1+\theta_3+\pi/2$ from the first condition, the second condition becomes
\begin{equation}\label{eq:n4tan}
\cos(\theta_1+\theta_3-2\theta_0)\cos(\theta_4-\theta_2)=\sin(\theta_1+\theta_3-2\theta_0)\cos(\theta_3-\theta_1).
\end{equation}
If $\cos(\theta_1+\theta_3-2\theta_0)\neq0$, this can also be written as $\cos(\theta_4-\theta_2)=\tan(\theta_1+\theta_3-2\theta_0)\cos(\theta_3-\theta_1)$.

\begin{example}[Mirror-symmetric cuts]\label{ex:mirror}
Take $\theta_2=\theta_0$, $\theta_4=\theta_0+\pi/2$ and $\theta_{1,3}=\theta_0\mp\delta$ with $0<\delta<\pi/2$. One chord passes through $C$, one is perpendicular to $OC$, and the other two are mirror images of each other in the line $OC$. Both conditions in~\eqref{eq:n4special} hold. Geometrically, reflection in the line $OC$ maps the arrangement to itself and exchanges odd and even pieces ($S_1\leftrightarrow S_2$, $S_3\leftrightarrow S_8$, and so on), so the equality is forced by symmetry. When $\delta\neq\pi/4$ the chords are \emph{not} equally spaced.
\end{example}

\subsection*{Universal fairness}

\begin{proposition}[Universal fairness, $n$ even]\label{prop:universal}
Let $n$ be even. A set of directions is universally fair if and only if
\[
A=\frac{\pi}{2}\qquad\text{and}\qquad Z=\sum_{j=1}^{n}(-1)^{j}e^{2i\theta_j}=0 .
\]
\end{proposition}

\begin{proof}
By Theorem~\ref{thm:even}, $\Sodd-\Seven=2a^{2}(A-\pi/2)+r_0^{2}\Imag(e^{-2i\theta_0}Z)$, and this must vanish for all $(r_0,\theta_0)\in[0,a)\times\mathbb R$. Setting $r_0=0$ forces $A=\pi/2$. Then $\Imag(e^{-2i\theta_0}Z)=0$ for all $\theta_0$, which forces $Z=0$. The converse is immediate.
\end{proof}

\begin{corollary}[Goldberg~\cite{Goldberg}]\label{cor:goldberg}
If $n\ge4$ is even and the chords are equally spaced, $\theta_j=\theta_1+(j-1)\pi/n$, then the directions are universally fair. For $n=4$ this is the pizza theorem of Upton~\cite{Upton}.
\end{corollary}

\begin{proof}
Here $A=\frac n2\cdot\frac{\pi}{n}=\frac{\pi}{2}$. With $\omega=e^{2\pi i/n}$,
$Z=-e^{2i\theta_1}\sum_{m=0}^{n-1}(-\omega)^{m}=-e^{2i\theta_1}\,\frac{1-(-\omega)^{n}}{1+\omega}=0$,
because $(-\omega)^{n}=1$ for even $n$ and $\omega\neq-1$ for $n\neq2$.
\end{proof}

For four chords the converse also holds.

\begin{proposition}[Converse of the pizza theorem for four chords]\label{prop:n4unique}
Four directions are universally fair if and only if they are equally spaced, $\theta_{j+1}-\theta_j=\pi/4$.
\end{proposition}

\begin{proof}
Sufficiency is Corollary~\ref{cor:goldberg}. For necessity, write $u_j=e^{2i\theta_j}$. These are four distinct points on the unit circle, because the lines are distinct. The condition $Z=0$ says $u_1+u_3=u_2+u_4=:s$. If $s\neq0$, the unit vectors $u$ with $|s-u|=1$ are exactly the two intersection points of the unit circles centred at $0$ and at $s$. Therefore $\{u_1,u_3\}=\{u_2,u_4\}$, which contradicts distinctness. Hence $s=0$, so $u_3=-u_1$ and $u_4=-u_2$, that is, $\theta_3=\theta_1+\pi/2$ and $\theta_4=\theta_2+\pi/2$. Then $A=2(\theta_2-\theta_1)=\pi/2$ gives $\theta_2-\theta_1=\pi/4$.
\end{proof}

\section{Four pieces ($n=2$)}\label{sec:four}

For two chords we distinguish two angles. The first is the \emph{oriented opening} $\Delta=\theta_2-\theta_1\in(0,\pi)$ of the piece $S_1$, measured counter-clockwise between consecutive rays. The second is the \emph{acute angle between the lines}, $\psi=\min(\Delta,\pi-\Delta)\in(0,\pi/2]$. Let $\sigma=\varphi_1+\varphi_2$; then $\sigma/2$ is the angle from $OC$ to the bisector of $S_1$. Theorem~\ref{thm:even} gives
\begin{equation}\label{eq:n2}
S_1+S_3=a^{2}\Delta+r_0^{2}\sin\Delta\cos\sigma,\qquad S_2+S_4=a^{2}(\pi-\Delta)-r_0^{2}\sin\Delta\cos\sigma,
\end{equation}
and equality holds if and only if
\begin{equation}\label{eq:n2cond}
a^{2}\Bigl(\Delta-\frac{\pi}{2}\Bigr)+r_0^{2}\sin\Delta\cos\sigma=0
\quad\Longleftrightarrow\quad
\frac{r_0^{2}}{2}\bigl[\sin2\varphi_2-\sin2\varphi_1\bigr]+a^{2}\Bigl(\Delta-\frac{\pi}{2}\Bigr)=0 .
\end{equation}

\begin{corollary}[Perpendicular chords]\label{cor:perp}
Two perpendicular chords through $O\neq C$ give equal alternate sums if and only if one of them passes through the centre $C$.
\end{corollary}

\begin{proof}
With $\Delta=\pi/2$, condition~\eqref{eq:n2cond} becomes $\sin2\varphi_2=\sin2\varphi_1$. Since $\varphi_2=\varphi_1+\pi/2$, this reads $-\sin2\varphi_1=\sin2\varphi_1$, so $\sin2\varphi_1=0$. Hence $\varphi_1\in\{0,\pi/2\}\pmod\pi$, which means one of the two chords lies along $OC$.
\end{proof}

\begin{proposition}[Admissible angles for a fair two-chord cut]\label{prop:n2window}
Let $0<\kappa<1$, and let $u^{*}(\kappa)\in(0,\kappa^{2})$ be the unique root of $u=\kappa^{2}\cos u$ in $[0,\pi/2)$. Two chords through $O$ can be rotated together about $O$ into a position where the alternate sums are equal if and only if
\[
\Bigl|\Delta-\frac{\pi}{2}\Bigr|=\frac{\pi}{2}-\psi\;\le\; u^{*}(\kappa)\;<\;\Bigl(\frac{r_0}{a}\Bigr)^{2}.
\]
In words, the acute angle between the chords must lie within $u^*(\kappa)$ of a right angle. The criterion is the same whether it is expressed through $\Delta$ or through $\psi$. No two-chord direction set is universally fair.
\end{proposition}

\begin{proof}
By~\eqref{eq:n2cond}, a fair orientation exists if and only if $\cos\sigma=(\pi/2-\Delta)/(\kappa^{2}\sin\Delta)$ has a solution $\sigma$, that is, if and only if $|u|\le\kappa^{2}\cos u$ with $u=\pi/2-\Delta\in(-\pi/2,\pi/2)$. Since $\cos u$ is even and $u\mapsto u-\kappa^{2}\cos u$ is strictly increasing on $[0,\pi/2)$, this is equivalent to $|u|\le u^*$. We have $|u|=\pi/2-\psi$, and $u^*<\kappa^2$ because $\cos u^*<1$. Finally, $Z=e^{2i\theta_2}-e^{2i\theta_1}\neq0$, so Proposition~\ref{prop:universal} rules out universal fairness.
\end{proof}

\begin{remark}[The extreme configurations]\label{rem:n2endpoints}
Equality $|\Delta-\pi/2|=u^*$ requires $\cos\sigma=\pm1$, so the line $OC$ bisects one pair of opposite pieces. If $\cos\sigma=1$, then $\Delta=\pi/2-u^*$ and the bisector of $S_1$ points at $C$. If $\cos\sigma=-1$, then $\Delta=\pi/2+u^*$; the bisector of $S_1$ is perpendicular to $OC$, and the line $OC$ bisects $S_2$ and $S_4$. In both cases the configuration is symmetric in the line $OC$, and the \emph{acute} pair of pieces, with opening $\psi=\pi/2-u^*$, is the pair containing $C$. Their extra radial length exactly compensates for their smaller angle. Strictly inside the window there are two fair orientations, $\pm\sigma$, mirror images of each other in $OC$.
\end{remark}

\begin{figure}[t]
\centering
\includegraphics[width=\textwidth]{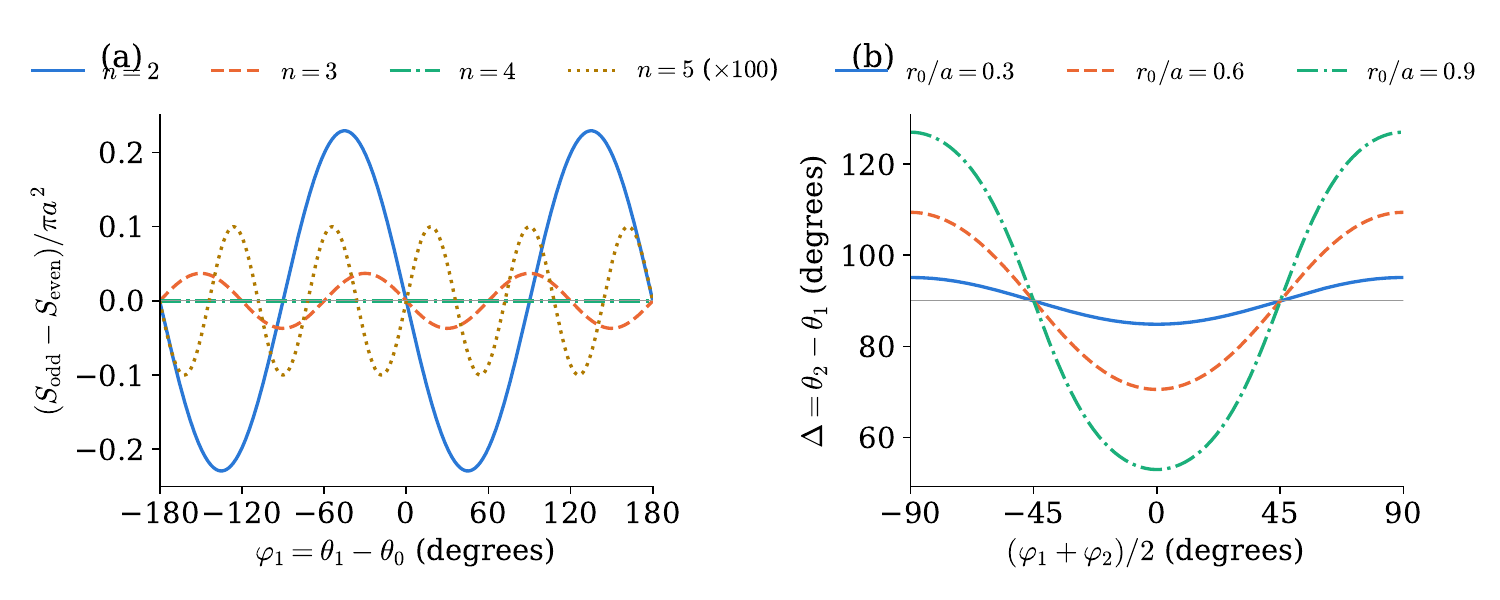}
\caption{(a)~Normalised imbalance $(\Sodd-\Seven)/\pi a^{2}$ for $n$ equally spaced chords as a function of the orientation $\varphi_1$, for $r_0/a=0.6$. For $n=5$ the plotted ordinate is $100\times(\Sodd-\Seven)/\pi a^{2}$; the other curves are unscaled. The $n=4$ curve vanishes identically (Corollary~\ref{cor:goldberg}). The $n=3$ and $n=5$ curves vanish when a chord passes through $C$ ($\sin n\varphi_1=0$) and have opposite signs near the centre (Proposition~\ref{prop:asym}). (b)~Fair two-chord cuts: the oriented opening $\Delta$ of $S_1$ against the direction $(\varphi_1+\varphi_2)/2$ of its bisector relative to $OC$, from~\eqref{eq:n2cond}. Both axes are unscaled and in degrees.}
\label{fig:imbalance}
\end{figure}

Figure~\ref{fig:imbalance}(b) shows the curve of fair cuts in the $(\sigma/2,\Delta)$ plane for several values of $r_0/a$. All curves pass through $(\pm\pi/4,\pi/2)$, which are the perpendicular cuts of Corollary~\ref{cor:perp}.

\section{Six pieces ($n=3$)}\label{sec:six}

For three chords ($\theta_{4}=\theta_1+\pi$, $\theta_5=\theta_2+\pi$, $\theta_6=\theta_3+\pi$), Theorem~\ref{thm:odd} gives
\begin{equation}\label{eq:n3}
S_1+S_3+S_5=\frac{\pi a^{2}}{2}-a^{2}\bigl[h(\kappa\sin\varphi_1)-h(\kappa\sin\varphi_2)+h(\kappa\sin\varphi_3)\bigr].
\end{equation}
In terms of $x_j=\arcsin(\kappa\sin\varphi_j)$ this reads
\[
S_1+S_3+S_5=\frac{\pi a^{2}}{2}+a^{2}\Bigl[(x_2-x_1-x_3)+\tfrac12(\sin2x_2-\sin2x_1-\sin2x_3)\Bigr].
\]
There is no $r_0^{2}\sin2\varphi$ term. The contributions of $S_1$, $S_3$ and $S_5$ to such terms are $\tfrac14r_0^2(\sin2\varphi_2-\sin2\varphi_1)$, $\tfrac14r_0^2(\sin2\varphi_1-\sin2\varphi_3)$ and $\tfrac14r_0^2(\sin2\varphi_3-\sin2\varphi_2)$, which sum to zero. Thus the six-piece sums are equal if and only if
\begin{equation}\label{eq:n3cond}
h(\kappa\sin\varphi_1)+h(\kappa\sin\varphi_3)=h(\kappa\sin\varphi_2).
\end{equation}

\begin{corollary}\label{cor:n3}
For three chords the following hold.
\begin{enumerate}
\item If $O=C$, the alternate sums are equal for all angles $\theta_1,\theta_2,\theta_3$.
\item If $O\neq C$ and one chord passes through $C$, the sums are equal if and only if the other two chords are mirror images of each other in that chord.
\item To first order in $\kappa$, condition~\eqref{eq:n3cond} reads $\sin\varphi_1-\sin\varphi_2+\sin\varphi_3=0$, that is, $d_1-d_2+d_3=0$ for the offsets of Lemma~\ref{lem:half}.
\end{enumerate}
\end{corollary}

\begin{proof}
Part (i) is Remark~\ref{rem:centred}. For (ii), suppose first that $\varphi_2=0$. Since $h$ is odd and strictly increasing, \eqref{eq:n3cond} gives $\sin\varphi_3=-\sin\varphi_1$, so $\varphi_3\equiv-\varphi_1\pmod\pi$; the other solution would make two chords coincide. If instead $\varphi_1=0$, we obtain $\sin\varphi_3=\sin\varphi_2$, so $\varphi_3\equiv-\varphi_2\pmod\pi$. The case $\varphi_3=0$ is analogous. Part (iii) uses $h(s)=2s+O(s^{3})$.
\end{proof}

\section{Odd numbers of chords: imbalance near the centre and impossibility of universal fairness}\label{sec:asym}

For odd $n$ and equally spaced chords, $\varphi_j=\varphi_1+(j-1)\pi/n$, the low-order terms of~\eqref{eq:odd} cancel exactly. This explains the small amplitudes in Figure~\ref{fig:imbalance}(a). Gluck~\cite{Gluck} has given the complete Fourier series of the imbalance, together with a global bound. The following proposition extracts its leading term directly from Theorem~\ref{thm:odd}, with an explicit remainder bound that is uniform in the orientation.

\begin{proposition}\label{prop:asym}
Let $n\ge1$ be odd and the chords equally spaced. Put $m_0=(n-1)/2$ and
\[
L_n(\kappa,\varphi_1)=-2^{\,3-n}\binom{1/2}{m_0}\,a^{2}\kappa^{n}\sin(n\varphi_1).
\]
Then, for all $\varphi_1\in\mathbb R$ and $0\le\kappa<1$,
\begin{equation}\label{eq:asym}
\bigl|\Sodd-\Seven-L_n(\kappa,\varphi_1)\bigr|\le\frac{4n\,a^{2}\kappa^{n+2}}{1-\kappa^{2}}.
\end{equation}
For example, $L_3=-\tfrac12a^2\kappa^3\sin3\varphi_1$.
\end{proposition}

\begin{proof}
Write $h(s)=\sum_{m\ge0}c_ms^{2m+1}$ with $c_m=2(-1)^{m}\binom{1/2}{m}/(2m+1)$. This series converges for $|s|<1$ because $h(s)=2\int_0^s\sqrt{1-t^2}\,dt$, and $|c_m|\le2$. Let $P_m(\varphi_1)=\sum_{j=1}^n(-1)^{j+1}\sin^{2m+1}\varphi_j$, so that $\sum_j(-1)^{j+1}h(\kappa\sin\varphi_j)=\sum_mc_m\kappa^{2m+1}P_m$. The power $\sin^{2m+1}\varphi$ is a linear combination of $\sin(q\varphi)$ with $q$ odd and $q\le2m+1$. The coefficient of the top harmonic $\sin((2m+1)\varphi)$ is $(-1)^m2^{-2m}$. For odd $q$,
\[
\sum_{j=1}^{n}(-1)^{j+1}\sin\bigl(q\varphi_1+q(j-1)\pi/n\bigr)=\Imag\Bigl(e^{iq\varphi_1}\sum_{l=0}^{n-1}\rho^{l}\Bigr),\qquad\rho=-e^{iq\pi/n}.
\]
Since $\rho^{n}=1$, this vanishes unless $\rho=1$, that is, unless $q$ is an odd multiple of $n$, in which case it equals $n\sin(q\varphi_1)$. Hence $P_m\equiv0$ for $m<m_0$, and $P_{m_0}=(-1)^{m_0}2^{1-n}n\sin(n\varphi_1)$ exactly. For $m>m_0$ we use only $|P_m|\le n$. Therefore
\[
\Bigl|\sum_j(-1)^{j+1}h(\kappa\sin\varphi_j)-c_{m_0}\kappa^nP_{m_0}\Bigr|\le\sum_{m>m_0}2n\kappa^{2m+1}=\frac{2n\kappa^{n+2}}{1-\kappa^2}.
\]
Multiply by $-2a^2$, use~\eqref{eq:odd}, and note that $-2a^2c_{m_0}\kappa^nP_{m_0}=L_n$ because $(-1)^{m_0}c_{m_0}n=2\binom{1/2}{m_0}$.
\end{proof}

\begin{remark}[Scope of the sign statement]\label{rem:sign}
The centre $C$ lies in $S_1$ exactly when $-\pi/n<\varphi_1<0$, and then $\sin n\varphi_1<0$. We have $\binom{1/2}{m}>0$ for $m=0$ and for odd $m$, and $\binom{1/2}{m}<0$ for even $m\ge2$. By~\eqref{eq:asym}, $\Sodd-\Seven$ has the sign of $L_n$ whenever
\[
2^{3-n}\Bigl|\binom{1/2}{m_0}\Bigr|\,|\sin n\varphi_1|>\frac{4n\,\kappa^{2}}{1-\kappa^{2}}.
\]
On this set, the share containing $C$ is larger when $n=1$ or $n\equiv3\pmod4$, and smaller when $n\equiv1\pmod4$ with $n\ge5$. This conclusion is \emph{not uniform} in the orientation: as $\varphi_1$ approaches a zero of $\sin n\varphi_1$ (a chord through $C$, where the two shares are exactly equal by symmetry), the admissible range of $\kappa$ shrinks to zero. Proposition~\ref{prop:asym} therefore explains the leading-order mechanism but does \emph{not} prove the global sign rule. That rule, valid for all interior points not on a chord through $C$, is the theorem of Mabry and Deiermann~\cite{MabryDeiermann}.
\end{remark}

Finally, no direction set with an odd number of chords is universally fair, whether or not the chords are equally spaced.

\begin{proposition}\label{prop:oddnever}
If $n$ is odd, no set of directions $\theta_1<\dots<\theta_n<\theta_1+\pi$ is universally fair.
\end{proposition}

\begin{proof}
Suppose the directions are universally fair. By Theorem~\ref{thm:odd},
\[
G(\kappa,\theta_0)=\sum_{j=1}^n(-1)^{j+1}h\bigl(\kappa\sin(\theta_j-\theta_0)\bigr)=\sum_{m\ge0}c_m\kappa^{2m+1}P_m(\theta_0),\qquad P_m(\theta_0)=\sum_{j=1}^n(-1)^{j+1}\sin^{2m+1}(\theta_j-\theta_0),
\]
vanishes for all $\kappa\in[0,1)$ and all $\theta_0$. For fixed $\theta_0$ the right-hand side is a power series in $\kappa$ that converges on $(-1,1)$ and vanishes on $[0,1)$. By the identity theorem for power series, each coefficient vanishes, and since $c_m\neq0$ we get $P_m\equiv0$ for every $m\ge0$.

Next we show by induction that $T_q(\theta_0)=\sum_j(-1)^{j+1}\sin(q(\theta_j-\theta_0))\equiv0$ for every odd $q$. For $m=0$, $P_0=T_1$. Suppose $T_1\equiv T_3\equiv\dots\equiv T_{2m-1}\equiv0$. Write $\sin^{2m+1}x=\sum_{l=0}^{m}\beta_{m,l}\sin((2l+1)x)$ with $\beta_{m,m}=(-1)^m2^{-2m}\neq0$. Then $0\equiv P_m=\sum_{l\le m}\beta_{m,l}T_{2l+1}=\beta_{m,m}T_{2m+1}$, so $T_{2m+1}\equiv0$. Since $T_q(\theta_0)=\Imag\bigl(e^{-iq\theta_0}W_q\bigr)$ with $W_q=\sum_j(-1)^{j+1}e^{iq\theta_j}$, vanishing for all $\theta_0$ means $W_q=0$ for all odd $q$.

Finally, put $z_j=e^{i\theta_j}$ and $y_j=z_j^2=e^{2i\theta_j}$. The $y_j$ are pairwise distinct because the lines are distinct. Taking $q=2k+1$ with $k=0,1,\dots,n-1$ gives the linear system $\sum_{j=1}^n y_j^{\,k}\,\bigl[(-1)^{j+1}z_j\bigr]=0$, $0\le k\le n-1$. Its matrix is the Vandermonde matrix $(y_j^{\,k})$, which is invertible because the $y_j$ are distinct. Hence $(-1)^{j+1}z_j=0$ for all $j$, which is absurd since $|z_j|=1$.
\end{proof}

\section{Conclusions}\label{sec:conclusions}

For $n$ concurrent chords with arbitrary directions, the alternate sums of sector areas are equal exactly when~\eqref{eq:evencond} holds ($n$ even) or~\eqref{eq:oddcond} holds ($n$ odd). For the cases treated in detail:
\begin{align*}
n=4\ (\text{8 pieces}):&\quad \tfrac12r_0^{2}\bigl[\sin2\varphi_2-\sin2\varphi_1+\sin2\varphi_4-\sin2\varphi_3\bigr]+a^{2}\bigl(\theta_2-\theta_1+\theta_4-\theta_3-\tfrac{\pi}{2}\bigr)=0,\\
n=2\ (\text{4 pieces}):&\quad \tfrac12r_0^{2}\bigl[\sin2\varphi_2-\sin2\varphi_1\bigr]+a^{2}\bigl(\theta_2-\theta_1-\tfrac{\pi}{2}\bigr)=0,\\
n=3\ (\text{6 pieces}):&\quad h(\kappa\sin\varphi_1)-h(\kappa\sin\varphi_2)+h(\kappa\sin\varphi_3)=0 .
\end{align*}
When $O$ coincides with $C$, these reduce to $(\theta_2-\theta_1)+(\theta_4-\theta_3)=\pi/2$, to $\theta_2-\theta_1=\pi/2$, and to no condition at all, respectively.

The difference between the even and odd cases has a simple cause. For even $n$ the alternate sums consist of opposite pairs, whose areas follow from the power of a point; the condition is then affine in $r_0^2$. For odd $n$ they are alternating sums of half-discs, whose areas depend nonlinearly on the distances of the chords from the centre. These general theorems contain the classical pizza theorem as a special case. They also give new statements for unequal angles: fairness along a diameter, the rigidity of the equally spaced four-chord cut, the impossibility of universally fair odd cuts, and the narrow window of fair two-chord cuts. Natural next steps include the corresponding statements for the crust (perimeter) and for toppings~\cite{MabryDeiermann,Hirschhorn}. Another is the classification of universally fair direction sets for even $n\ge6$, where $A=\pi/2$ and $Z=0$ admit solutions other than equal spacing.

\appendix
\section{Reproducible numerical verification}\label{app:numerics}

The principal closed-form area identities (Proposition~\ref{prop:sector}, Lemmas~\ref{lem:pair} and~\ref{lem:half}, Theorems~\ref{thm:even} and~\ref{thm:odd}) and selected geometric special cases were checked independently against direct adaptive quadrature of $\tfrac12\int r^2(\theta)\,d\theta$, with $r$ given by~\eqref{eq:r}. The universally quantified statements, such as Propositions~\ref{prop:universal}, \ref{prop:n4unique} and~\ref{prop:oddnever}, rest on their proofs alone and are not covered by sampling. The checks corroborate the algebra but do not replace the proofs. The script \texttt{verify\_all.py}, supplied with the article, runs all tests in about two seconds. From the root of the supplementary package:
\begin{center}
\texttt{python3 scripts/verify\_all.py}\qquad(defaults: \texttt{--N 500 --seed 20261027}; requires NumPy and SciPy).
\end{center}

\emph{Sampling protocol.} The radius was drawn as $a\sim U[0.5,2]$ and the ratio as $\kappa\sim U[0,0.95]$; in addition, every tenth sample uses $\kappa=0$ exactly. The direction of the centre was drawn as $\theta_0\sim U[-\pi,\pi)$. The chord directions were $n$ independent $U[0,\pi)$ variates, sorted. A sample was rejected and redrawn unless every consecutive gap, including the wrap-around gap $\theta_1+\pi-\theta_n$, was at least $0.01$ rad, so the order is strict. Quadrature used \texttt{scipy.integrate.quad} with $\texttt{epsabs}=\texttt{epsrel}=10^{-13}$. Errors are reported relative to the disc area, as $\max|\text{closed form}-\text{quadrature}|/\pi a^2$, and a test passes if this is below $10^{-11}$.

\begin{table}[h]
\centering
\small
\begin{tabular}{llrc}
\toprule
Test & Result checked & Samples & Max.\ relative error\\
\midrule
T1 & Prop.~\ref{prop:sector}, arbitrary $\alpha<\beta$ & 500 & $4.2\times10^{-16}$\\
T2 & Lemma~\ref{lem:pair}, $n\in\{1,\dots,6\}$ & 500 & $6.1\times10^{-16}$\\
T3 & Lemma~\ref{lem:half}, $n\in\{1,\dots,6\}$ & 500 & $4.4\times10^{-16}$\\
T4 & Thm.~\ref{thm:even}, $n\in\{2,4,6,8\}$ & 500 & $4.8\times10^{-16}$\\
T5 & Thm.~\ref{thm:odd}, $n\in\{1,3,5,7\}$ & 500 & $6.0\times10^{-16}$\\
T6 & Example~\ref{ex:mirror} & 100 & $7.4\times10^{-16}$\\
T7 & Cor.~\ref{cor:goldberg}, $n\in\{4,6,8\}$ & 100 & $1.1\times10^{-15}$\\
T8 & Cor.~\ref{cor:perp} & 100 & $5.0\times10^{-16}$\\
T9 & Prop.~\ref{prop:n2window}, both window endpoints & 200 & $2.8\times10^{-13}$\\
T10 & Cor.~\ref{cor:n3}(i)--(ii) & 300 & $7.6\times10^{-16}$\\
\bottomrule
\end{tabular}
\caption{Numerical checks produced by \texttt{verify\_all.py} with the default settings.}
\label{tab:numerics}
\end{table}

The larger error in T9 comes from the root-finding tolerance used to compute $u^*(\kappa)$, not from the quadrature. Test T11 checks Proposition~\ref{prop:asym} for $n\in\{1,3,5,7\}$ and $\kappa\in\{0.05,0.1,0.2,0.4\}$ on a grid of $721$ orientations. The ratio of the exact imbalance to $L_n$ tends to $1$ as $\kappa\to0$; for example, at $\kappa=0.05$ it is $1.0005$ for $n=3$ and $1.0016$ for $n=5$. The largest observed value of $|\Sodd-\Seven-L_n|$ divided by the bound in~\eqref{eq:asym} was $0.17$, so the bound holds with room to spare.

\end{document}